\documentclass{amsproc}
\usepackage{amssymb}

\usepackage{graphicx}
\usepackage{amscd}

\theoremstyle{plain}
\newtheorem{acknowledgement}{Acknowledgement}
\newtheorem{corollary}{Corollary}

\newtheorem{proposition}{Proposition}
\newtheorem{remark}{Remark}
\newtheorem{theorem}{Theorem}
\numberwithin{equation}{section}
\newtheorem{definition}{Definition}

\begin{document}
\title[Banach Spaces of Large Density Character]{On Banach Spaces of Large Density Character and Unconditional Sequences}
\author{Eugene Tokarev}
\address{B.E. Ukrecolan, 33-81, Iskrinskaya str.,61005, Kharkiv-5, Ukraine}
\email{evtokarev@yahoo.com}
\dedicatory{Dedicated to the memory of S. Banach.}

\begin{abstract}
It is shown that any Banach space $X$ of density $dens(X)\geq \gimel
_{\omega _{1}}$ contains an (infinite) unconditional sequence. Here for any
ordinal $\alpha $ the cardinal $\gimel _{\alpha }$ is given by induction: $%
\gimel _{0}=\aleph _{0}$; $\gimel _{\alpha +1}=\exp (\gimel _{\alpha })$; $%
\gimel _{\alpha }=\sup \{\gimel _{\beta }:\beta <\alpha \}$ if $\alpha $ is
a limit ordinal, and $\omega _{1}$ is the first uncountable ordinal.

Every Banach space of density $dens(X)=\tau \geq \omega _{1}$ contains an
unconditional sequence of cardinality $\tau $ if and only if $\tau $ is a
weakly compact cardinal.

Other results that are concerned large cardinals and unconditional sequences
are also presented.
\end{abstract}

\maketitle

\section{Introduction.}

The known J. Lindenstrauss' problem [1] was:

\textit{Whether every infinite dimensional Banach space contains an
(infinite) unconditional sequence?}

It was successfully solved (in negative) by W.T. Gowers and B. Maurey [2].

More formerly, J. Ketonen [3] has proved that any Banach space $X$ of
dimension\ $dim(X)=\varrho $, where $\varrho $ is a Ramsey cardinal
(definitions see below) contains a subspace with an unconditional basis of
the cardinality $\varrho $.

Recall that a \textit{dimension} $dim(X)$ of a Banach space $X$ is the least
cardinality of a subset $A\subset X$, which \textit{linear span} $lin(A)$ is
dense in $X$ (equivalently: such $A\subset X$ that a \textit{closure} $%
\overline{lin(A)}$, which will be in the future denoted $span(A)$, is the
whole space $X$). If $X$ is of infinite dimension then its dimension $dim(X)$
is exactly equal to its \textit{density character} $dens(X)$ - the least
cardinality of a subset $B\subset X$, which is dense in $X$.

The result of [3] was improved in [4] (see also [5]), where was shown that
for any cardinal $\tau $ there exists a cardinal $\varkappa (\tau )$
(namely, an P. Erd\"{o}s' cardinal $\varkappa \left( \tau \right)
\rightarrow \left( \tau \right) _{2}^{<\omega }$) such that any Banach space
of dimension $\geq \varkappa \left( \tau \right) $ contains an unconditional
sequence of cardinality $\tau $. Cardinals of kind $\varkappa (\tau )$ are
too large. It is known that $\iota <\varkappa \left( \omega \right)
<\varkappa \left( \omega _{1}\right) <...<\varkappa \left( \varkappa \right)
=\varrho $, where $\iota $ is the first inaccessible cardinal, and $\varrho $
is the first Ramsey cardinal (here and below $\omega $ is the first infinite
cardinal; $\omega _{1}$ is the first uncountable cardinal). Moreover, each $%
\varkappa (\tau )$ is inaccessible itself. So, the existence of $\varkappa
\left( \omega \right) $, of $\varkappa \left( \omega _{1}\right) $, etc.
must be postulated separately.

P. Terenzi [6] asked:

\textit{Whether there exists a cardinal }$\tau $ (\textit{without any
additional set-theoretical assumptions) such that any Banach space \ }$X$
\textit{of dimension }$dim(X)\geq \tau $ \textit{contains an infinite
unconditional sequence?}

By using of some methods from the model theory it may be proved that the
existence of $\varkappa \left( \omega \right) $ implies that any Banach
space of dimension $\geq \gimel _{\omega _{1}}$ contains an unconditional
sequence (the cardinal $\gimel _{\omega _{1}}$, which existence does not
depend on any set-theoretical hypotheses, will be defined later).

Below it will be shown that, independently of the existence (or
non-existence) of $\varkappa \left( \omega \right) $, any Banach space $X$
of $dim(X)\geq \gimel _{\omega _{1}}$ contains an unconditional sequence.

A sequence $\{x_{\alpha }:\alpha <\tau \}$ of elements of a Banach space $X$
will be called a\textit{\ }$\tau $\textit{-sequence}$.$ The following
question is of a separate interest:

\textit{What cardinals} $\tau $ \textit{has the property: each Banach space}
$X$ \textit{of dimension} $\tau $ \textit{contains an unconditional} $\tau $%
\textit{-sequence?}

It will be shown that a cardinal $\tau $ has a connected property (namely,
each Banach space of dimension $\tau $ contains an subsymmetric $\tau $%
-sequence) if and only if $\tau $ is a\ weakly compact cardinal. Hence, the
same property of Ramsey cardinals (cf. [1]) is a consequence of their weak
compactness.

By the way, some partial results in the S. Banach's problem on the separable
quotient will be obtained.

\section{Definitions and notations}

\textit{Ordinals} will be denoted by small Greece letters $\alpha ,\beta
,\gamma $. \textit{Cardinals} are identified with the least ordinals of a
given cardinality and are denoted either $\iota ,\tau ,\varkappa ,\varrho
,\sigma $, or by using Hebrew letters $\aleph $, $\gimel $ (may be with
indices). As usual, $\omega $ and $\omega _{1}$ denote respectively the
first infinite and the first uncountable cardinals (=ordinals).

For a cardinal $\tau $ its \textit{predecessor} (i.e. the least cardinal
that is strongly greater then $\tau $) is denoted by $\tau ^{+}$.

\textit{The} \textit{confinality of} $\tau ,$ $cf(\tau )$ is the least
cardinality of a set $A\subset \tau $ such that $\tau =\sup A$.

Let $A$, $B$ be sets. A symbol $^{B}A$ denotes a set of all functions from $%
B $ to $A$.

In a general case, a cardinality of the set $^{B}A$ is denoted either $%
card(A)^{card(B)}$ or $\varkappa ^{\tau }$, if $card(A)=\varkappa $; $%
card(B)=\tau $.

A symbol $\exp (\tau )$ (or, equivalently, $2^{\tau }$) denotes a
cardinality of a set $^{\tau }2$ of all subsets of $\tau $.

Let $J$ be a set; $n<\omega $. Let
\begin{equation*}
\lbrack J]^{n}=\{A\subset J:card(A)=n\}
\end{equation*}

$J$ may be considered as a well-ordered set (say, by a relation '$<$'). So, $%
[J]^{n}$ may be regarded as a set of all $n$-tuples $a_{0}<a_{1}<...<a_{n-1}$
where $a_{i}\in J$ ($i<n$).

For a cardinal $\sigma $ a function $P:[J]^{n}\rightarrow \sigma $ is called
a\textit{\ partition} of $[J]^{n}$ (on $\sigma $ parts).

If $P:[J]^{n}\rightarrow \sigma $; $H\subset J$ and the function $P$ is a
constant on $[H]^{n}$, i.e., if all elements of $[H]^{n}$ belong to the same
class of the partition $P$, then the set $H$ is said to be \textit{%
homogeneous for} $P$.

The P. Erd\"{o}s' denotation
\begin{equation*}
\varkappa \rightarrow \left( \tau \right) _{\sigma }^{n}
\end{equation*}
will be used as a shortness for the following assertion:

\textit{For every partition} $P:[\varkappa ]^{n}\rightarrow \sigma $ \textit{%
there exists a homogeneous for }$P$ \textit{subset} $H\subset \varkappa $
\textit{of cardinality }$\tau $.

The another Erd\"{o}s' denotation
\begin{equation*}
\varkappa \rightarrow \left( \tau \right) _{\sigma }^{<\omega }
\end{equation*}
means:

\textit{For any family }$\{P_{n}:n<\omega \}$ \textit{of partitions} $%
P_{n}:[\varkappa ]^{n}\rightarrow \sigma $\textit{\ there exists a set} $%
H\subset \varkappa $, \textit{which is homogeneous for all of} $P_{n}$'s.

Let $\tau $ be a cardinal. An \textit{Erd\"{o}s' cardinal }$\varkappa \left(
\tau \right) $ is the least cardinal $\varkappa $ such that
\begin{equation*}
\varkappa \rightarrow \left( \tau \right) _{2}^{<\omega }.
\end{equation*}

A cardinal $\sigma $ is said to be \textit{weakly compact} if
\begin{equation*}
\sigma \rightarrow \left( \sigma \right) _{2}^{2}.
\end{equation*}

A cardinal $\rho $ is said to be a \textit{Ramsey cardinal }if
\begin{equation*}
\varrho \rightarrow \left( \varrho \right) _{2}^{<\omega }.
\end{equation*}

A cardinal $\upsilon $ is said to be \textit{measurable} if there exists a
countably additive function (measure) $\mu :\upsilon \rightarrow \{0,1\}$,
defined on all subsets of $\upsilon $, such that
\begin{equation*}
\mu \left( \upsilon \right) =1;\text{ }\mu \left( A\right) =0\text{ for
every subset }A\subset \upsilon \text{ of cardinality }card(A)<\upsilon .
\end{equation*}
Any Ramsey cardinal is weakly compact. Any measurable cardinal is weakly
compact and is a Ramsey cardinal.

Let $\alpha $ be an ordinal; $\tau $ be a cardinal. A cardinal $\gimel
_{\alpha }\left( \tau \right) $ is given by the induction:
\begin{equation*}
\gimel _{0}\left( \tau \right) =\tau ;\text{ \ \ }\gimel _{\alpha +1}\left(
\tau \right) =\exp \left( \gimel _{\alpha }\left( \tau \right) \right) ;
\end{equation*}
\begin{equation*}
\gimel _{\alpha }\left( \tau \right) =\sup \{\gimel _{\beta }\left( \tau
\right) :\beta <\alpha \},
\end{equation*}
if $\alpha $ is a\textit{\ limit ordinal} (i.e. if $\alpha $ is not of kind $%
\alpha =\beta +1$ for some ordinal $\beta $).

In a case $\tau =\omega $ ($=\aleph _{0})$ it will be written $\gimel
_{\alpha }$ instead of $\gimel _{\alpha }\left( \omega \right) $.

Recall that a cardinal $\aleph _{\alpha }$ is given by

\begin{center}
$\aleph _{0}=\omega $; $\aleph _{\alpha +1}=\left( \aleph _{\alpha }\right)
^{+}$; $\aleph _{\alpha }=\sup \{\aleph _{\beta }:\beta <\alpha \}$, if $%
\alpha $ is a limit ordinal.
\end{center}

The \textit{general continuum hypothesis} (GCH) is a conjecture:

\begin{center}
$\aleph _{\alpha }=\gimel _{\alpha }$ for every ordinal $\alpha $.
\end{center}

The \textit{continuum hypotheses} (CH) asserts that $\aleph _{1}=\gimel _{1}$%
.

In the future the following Erd\"{o}s-P. Rado theorem [7] will be used:

\textit{For every infinite cardinal} $\varkappa $
\begin{equation*}
\left( \gimel _{n}\left( \varkappa \right) \right) ^{+}\rightarrow \left(
\varkappa ^{+}\right) _{\varkappa }^{n+1}.
\end{equation*}

Let $\mathcal{B}$ be a (proper) class of all Banach spaces. Let $X\in
\mathcal{B}$; $dim(X)=\varkappa $; $\alpha $ be an infinite limit ordinal\
(so, $\alpha $ may be considered as a cardinal); $\alpha \leq \varkappa $.
An $\alpha $-sequence $\{x_{\beta }:\beta <\alpha \}$ of elements of $X$ is
said to be

\begin{itemize}
\item  \textit{Almost spreading,} if for any $n<\omega $, any $\varepsilon
>0 $, any scalars $\{a_{k}:k<n\}$ and any choosing of $%
i_{0}<i_{1}<...<i_{n-1}<\alpha $; $j_{0}<j_{1}<...<j_{n-1}<\alpha $ the
following inequality holds:
\begin{equation*}
(1-\varepsilon )\left\| \sum\nolimits_{k<n}a_{k}x_{j_{k}}\right\| \leq
\left\| \sum\nolimits_{k<n}a_{k}x_{i_{k}}\right\| \leq (1+\varepsilon
)\left\| \sum\nolimits_{k<n}a_{k}x_{j_{k}}\right\| ;
\end{equation*}

\item  \textit{Spreading, }if the previous is remain valid for $\varepsilon
=0$;

\item  $C$-\textit{unconditional}, where $C<\infty $ is a constant, if
\begin{equation*}
C^{-1}\left\| \sum\nolimits_{k<n}a_{k}\epsilon _{k}x_{i_{k}}\right\| \leq
\left\| \sum\nolimits_{k<n}a_{k}x_{i_{k}}\right\| \leq C\left\|
\sum\nolimits_{k<n}a_{k}\epsilon _{k}x_{i_{k}}\right\|
\end{equation*}
for any choosing of $n<\omega $; $\{a_{k}:k<n\}$; $i_{0}<i_{1}<...<i_{n-1}<%
\alpha $ and of signs $\{\epsilon _{k}\in \{+,-\}:k<n\}$.

\item  \textit{Unconditional}, if it is $C$-unconditional for some $C<\infty
$.

\item  \textit{Symmetric,} if for any $n<\omega $, any finite subset $%
I\subset \mathbb{\alpha }$ of cardinality $n$, any rearrangement $\varsigma $
of elements of $I$ and any scalars $\{a_{i}:i\in I\}$,
\begin{equation*}
\left\| \sum\nolimits_{i\in I}a_{i}z_{i}\right\| =\left\|
\sum\nolimits_{i\in I}a_{\varsigma (i)}z_{i}\right\| .
\end{equation*}

\item  \textit{Subsymmetric}, if it is both spreading and 1-unconditional.
\end{itemize}

Let $C<\infty $ be a constant. Two $\alpha $-sequences $\{x_{\beta }:\beta
<\alpha \}$ and $\{y_{\beta }:\beta <\alpha \}$ are said to be $C$\textit{%
-equivalent} if for any finite subset $I=\{i_{0}<i_{1}<...<i_{n-1}\}$ of $%
\alpha $ and for any choosing of scalars $\{a_{k}:k<n\}$
\begin{equation*}
C^{-1}\left\| \sum\nolimits_{k<n}a_{k}x_{i_{k}}\right\| \leq \left\|
\sum\nolimits_{k<n}a_{k}y_{i_{k}}\right\| \leq C\left\|
\sum\nolimits_{k<n}a_{k}x_{i_{k}}\right\| .
\end{equation*}

Two $\alpha $-sequences $\{x_{\beta }:\beta <\alpha \}$ and $\{y_{\beta
}:\beta <\alpha \}$ are said to be equivalent if they are $C$-equivalent for
some $C<\infty $.

It is known (cf. [8]) that for any spreading sequence $\{x_{n}:n<\omega \}$
the sequence of differences $\{x_{2m+1}-x_{2m}:m<\omega \}$ is
1-unconditional (and, hence, the same is true for any $\alpha $-sequence
when $2m+1$ and $2m$ are replaced with corresponding ordinals). Clearly, if $%
X\in \mathcal{B}$ contains an almost spreading $\alpha $-sequence then $X$
contains an\textit{\ almost unconditional }(i.e. $C$-unconditional for every
$C>1$) $\alpha $-sequence.

\begin{remark}
The same definitions may be used in a case when instead of $\alpha $%
-sequences will be regarded families $\{x_{i}:i\in I\}\subset X$ indexed by
elements of a linearly ordered set $\left\langle I,\ll \right\rangle $. In a
such case it will be said about spreading families, unconditional families
and so on.
\end{remark}

\section{On subspaces with unconditional bases in Banach spaces of large
dimension}

Let $X\in \mathcal{B}$; $dim(X)=\gimel _{\omega _{1}}$; $\{x_{\alpha
}:\alpha <\gimel _{\omega _{1}}\}$ be a sequence of linearly independent
elements of $X$ of norm $\left\| x_{\alpha }\right\| =1$ for all $\alpha $,
such that $X=span\{x_{\alpha }:\alpha <\gimel _{\omega _{1}}\}$. Let $%
span\{x_{\beta }:\beta <\gamma \}=X_{\gamma }$ for all $\gamma <$.$\gimel
_{\omega _{1}}$. It will be shown that $X$ contains an almost spreading
sequence (i.e., $\omega $-sequence). This would imply that $X$ contains an
unconditional sequence.

\begin{theorem}
Any Banach space of dimension $\gimel _{\omega _{1}}$ contains an infinite
dimensional subspace with an unconditional basis.
\end{theorem}

\begin{proof}
Let $dim(X)=\gimel _{\omega _{1}}$; $E_{\omega _{1}}=\{x_{\alpha }:\alpha
<\gimel _{\omega _{1}}\}$ - a sequence of linearly independent elements of $%
X $ of norm $\left\| x_{\alpha }\right\| =1$ for all $\alpha $ such that
\begin{equation*}
X=span\{x_{\alpha }:\alpha <\gimel _{\omega _{1}}\}
\end{equation*}
Let $\gamma <\omega _{1}$. Denote
\begin{equation*}
E_{\gamma }=\{x_{\alpha }:\alpha <\gimel _{\gamma }\};\text{ \ }X_{\gamma
}=span\{x_{\beta }:\beta <\gimel _{\gamma }\}.
\end{equation*}

Let $\Psi =\{\psi _{N}:N<\omega \}$ be a numeration of all finite sequences
of rational numbers (i.e., of elements of $\cup \{\mathbb{Q}^{n}:n<\omega \}
$).

For $\psi =\{a_{0},a_{1},...,a_{n}\}\in \Psi $ let $lh(\psi )=n=card\left(
\psi \right) -1$

It will be convenient to underline the correspondence between a number $N$
of a given $\psi _{N}\in \Psi $ in a mentioned numeration and its length $%
lh(\psi )$ in a following way.

Let
\begin{equation*}
\psi _{N}=\{a_{0}^{(N)},a_{1}^{(N)},...,a_{n_{N}}^{(N)},r^{(N)}\}.
\end{equation*}

To each $\psi _{N}\in \Psi $ will be assigned two formal symbols: $\phi
_{N}^{+}$ and $\phi _{N}^{-}$.

It will be said that a sequence $\frak{x}%
=\{x_{i_{0}},x_{i_{1}},...,x_{i_{N}}\}\subset E_{\omega _{1}}$, where $%
i_{0}<i_{1}<...<i_{N}$, satisfies $\phi _{N}^{+}$, symbolically: $\frak{x}%
\vDash \phi _{N}^{+}$ (resp., $\phi _{N}^{-}$; in symbols: $\frak{x}\vDash
\phi _{N}^{-}$) if
\begin{equation*}
\left\| \sum\nolimits_{k<n_{N}}a_{k}^{(N)}x_{i_{k}}\right\| \geq r^{(N)}
\end{equation*}
(respectively, if
\begin{equation*}
\left\| \sum\nolimits_{k<n_{N}}a_{k}^{(N)}x_{i_{k}}\right\| <r^{(N)}).
\end{equation*}
For any subset $G\subset E_{\omega _{1}}$ it will be written $G\vDash \phi
_{N}^{+}$ (resp., $G\vDash \phi _{N}^{-}$) if $\frak{x}\vDash \phi _{N}^{+}$
(resp., $\frak{x}\vDash \phi _{N}^{-}$) for all ordered $n_{N}$-tuples $%
\frak{x}=\{x_{i_{0}},x_{i_{1}},...,x_{i_{N}}\}$ from $G$.

Let us show by an induction that following is true:

\begin{itemize}
\item  \textit{There exists a sequence of signs} $\{\epsilon _{k}\in
\{+,-\}:k<n\}$\textit{, a confinal subset} $F_{n}\subset \omega _{1}$
\textit{and a sequence of subsets} $G_{\alpha }\subset E_{\alpha }$, $\alpha
\in F_{n}$ \textit{such that for any} $\alpha \in F_{n}$

\begin{enumerate}
\item  $card(G_{\alpha })\geq \gimel _{\gamma }$,\textit{\ where }$\alpha $
\textit{is }$\gamma $\textit{'s element of }$F_{n}$;

\item  $G_{\alpha }\vDash \phi _{m}^{\epsilon _{m}}$ \textit{for all} $m<n$.
\end{enumerate}
\end{itemize}

This assertion is obviously valid for $n=0$, $F_{1}=\omega _{1}$ and $%
G_{\alpha }=E_{\alpha }$.

Assume that $\epsilon _{0},\epsilon _{1},...,\epsilon _{n-1}$, $F_{n}$, $%
G_{\alpha }$ and $\alpha \in $ $F_{n}$ are already chosen. Let $\psi _{n}$
be of length $lh(\psi _{n)}=p_{n}\geq 1$.

Let $F_{n}^{\prime }\subset F_{n}$ consists of all elements of $F_{n}$ that
are of kind $\beta +(p_{n}-1)$ for some ordinal $\beta $:
\begin{eqnarray*}
F_{n}^{\prime } &=&\{\alpha \in F_{n}:\text{ there exists }\beta \in F_{n}%
\text{ \ such that }\left( \alpha \backslash \beta \right) \cap F_{n}\text{
\ } \\
&&\text{ \ \ \ \ \ \ \ \ \ consists exactly of\ \ }p_{n}-1\text{ \ elements}%
\}.
\end{eqnarray*}

Clearly, $F_{n}^{\prime }$ is confinal in $\omega _{1}$. Let $\alpha $ $\ $%
be a $\gamma $'s element of $F_{n}^{\prime }$. Then $card(G_{\alpha
})>\gimel _{\gamma +p_{n}-1}$.

The given $\psi _{n}$ generates a partition of $[G_{\alpha }]^{p_{n}}$ in
two parts, $G_{+}$ and $G_{-}$ in dependence of $\left( x_{i_{k}}\right)
_{k<p_{n}}\vDash \phi _{n}^{\epsilon }$, where $\epsilon \in \lbrack 0,1]$.
Hence, according to the Erd\"{o}s-Rado theorem, there exists a subset $%
H_{\alpha }\subset G_{\alpha }$ , which is homogeneous for this partition
and has a cardinality $card(H_{\alpha })=\gimel _{\gamma }$. Clearly, this
is true for every $\alpha \in F_{n}$. Since there exists only two
possibilities for $\epsilon \in \{+,-\}$, there exists a confinal subset $%
F_{n+1}\subset F_{n}$ such that for all its elements $\gamma \in F_{n+1}$
all $p_{m}$-tuples $\{x_{i_{0}},x_{i_{1}},...,x_{i_{p_{m}}}\}$ ($m\leq n+1$)
satisfy the same condition of kind $\phi _{m}^{\epsilon _{m}}$.

So, we choose $\epsilon $, $F_{n+1}$ and $\{H_{\alpha }:\alpha \in F_{n+1}\}$%
. This complete the induction.

As a result, it will be chosen a sequence of signs $\{\epsilon _{k}\in
\{+,-\}:k<\omega \}$, which, together with the given order on $\cup \{%
\mathbb{Q}^{n}:n<\omega \}$ and the norm of $X$ defines in a unique way (up
to an almost isometry) a subspace of $X$ with an approximately spreading
basis $\left( x_{n}\right) $ (for which may be chosen any subset of $\cap
F_{n}$). Of course, this intersection is non empty because of $%
card(F_{n})=\omega _{1}$ for all $n<\omega $ and since $cf(\omega _{1})\neq
\omega $.
\end{proof}

\begin{remark}
In an analogous way it may be proved that any Banach space of dimension $%
\geq \gimel _{\tau ^{+}}$ contains an unconditional $\tau $-sequence.
\end{remark}

It is conceivable that any Banach space of inaccessible dimension contains
an unconditional sequence of inaccessible cardinality.

The validity of this assertion is not known. A weaker version will be stated
below. Namely, it will be shown that if $\sigma $ is the first weakly
compact cardinal then any Banach space of dimension $\geq \sigma $ contains
an unconditional $\sigma $-sequence.

\section{On separable quotient spaces of Banach spaces}

There is a long standing Banach's problem:

\textit{Whether every infinite dimensional Banach space} $X$ \textit{\ has a
separable quotient space?}

It easy to point out a class $\mathcal{K}$ of cardinals such that any Banach
space of dimension $\varkappa \in \mathcal{K}$ has a separable quotient.

Let $X\in \mathcal{B}$; $dim(X)=\varkappa $; $\{x_{\alpha }:\left\|
x_{\alpha }\right\| =1;\alpha <\varkappa \}$ - a $\varkappa $-sequence of
linearly independent elements of $X$; $span\{x_{\alpha }:\alpha <\varkappa
\}=X$. Let $span\{x_{\beta }:\beta <\gamma \}=X_{\gamma }$ for all $\gamma
<\varkappa $.

Certainly, $X_{\alpha }\hookrightarrow X_{\gamma }$ (recall that a symbol $%
Y\hookrightarrow Z$ denotes that $Y$ is a subspace of $Z$) for any $\alpha
<\beta <\varkappa $. Obviously, $X$ may be represented as the closure of the
union of the chain $X_{0}\hookrightarrow X_{1}\hookrightarrow
...\hookrightarrow X_{\alpha }\hookrightarrow X_{\gamma }\hookrightarrow ...$%
:
\begin{equation*}
X=\overline{\cup \{X_{\alpha }:\alpha <\varkappa \}}.
\end{equation*}

\begin{theorem}
If a cardinal $\varkappa $ is of countable confinality, $cf(\varkappa
)=\omega $, then every Banach space of dimension $\varkappa $ has a
separable quotient space.
\end{theorem}

\begin{proof}
Let $cf(\varkappa )=\omega $. Then there exists a countable sequence $%
\varkappa _{0}<\varkappa _{1}<...<\varkappa _{n}<...<\varkappa $ such that $%
\sup \{\varkappa _{n}:n<\omega \}=\varkappa .$ Hence $X$ may be represented
as the closure of the union of a countable chain\ $X_{\varkappa
_{0}}\hookrightarrow ...\hookrightarrow X_{\varkappa _{n}}\hookrightarrow $:
\begin{equation*}
X=\overline{\cup \{X_{\varkappa _{n}}:n<\omega \}}.
\end{equation*}

According to [9], any Banach space that has a such representation has also a
separable quotient space.
\end{proof}

\begin{theorem}
Any $B$-convex Banach space $W$ of dimension $\geq \gimel _{\omega _{1}}$
has a separable quotient space.
\end{theorem}

\begin{proof}
Let $dim(W)\geq \gimel _{\omega _{1}}.$Then its conjugate also is of
dimension $dim(W^{\ast })\geq \gimel _{\omega _{1}}$ and thus contains an
unconditional sequence, which spans in $W^{\ast }$ either a reflexive
subspace, or a subspace isomorphic to $l_{1}$ or a subspace which is
isomorphic to $c_{0}$. Let $Z\hookrightarrow W^{\ast }$ be reflexive. Every
its separable subspace $Z_{0}\hookrightarrow Z$ is reflexive too. Since $%
Z_{0}$ is weakly* closed, clearly, $W$ has a quotient that is isometric to $%
\left( Z_{0}\right) ^{\ast }$.

If $W^{\ast }$ contains a subspace isomorphic to $c_{0}$ then (since $%
W^{\ast }$ is a conjugate space) $W^{\ast }$ contains a subspace isomorphic
to $L_{1}[0,1]$, which contains reflexive subspaces. Let $W^{\ast }$
contains a subspace $Y$ isomorphic to $l_{1}$. Consider its weak* closure $%
\left( \left( Y\right) _{\perp }\right) ^{\perp }$, where for a subset $%
A\subset W^{\ast }$ a symbol $\left( A\right) _{\perp }$ denotes its lower
annihilator
\begin{equation*}
\left( A\right) _{\perp }=\{w\in W:a(w)=0\text{ \ for all \ }a\in A\};
\end{equation*}
and for a subset $B\subset W$ a symbol $\left( B\right) ^{\perp }$ denotes
its upper annihilator
\begin{equation*}
\left( B\right) ^{_{\perp }}=\{w^{\ast }\in W^{\ast }:w^{\ast }(b)=0\text{ \
for all \ }b\in B\}.
\end{equation*}

Clearly, $\left( \left( Y\right) _{\perp }\right) ^{\perp }$ is an $%
L_{1}\left( \mu \right) $-space. If a measure $\mu $ is discrete then $%
W^{\ast }$ contains a weakly* closed subspace, isomorphic to $l_{1}$ and, of
course, $W$ has a separable quotient. If a measure $\mu $ contains a
nontrivial dispersed part, then $W^{\ast }$ contains a subspace isomorphic
to $L_{1}[0,1]$ and, as was mentioned before, has a separable quotient too.
\end{proof}

\section{Banach spaces of weakly compact dimension}

Recall that if a cardinal $\varkappa $ is weakly compact, i.e. if $\varkappa
\rightarrow \left( \varkappa \right) _{2}^{2}$, then $\varkappa \rightarrow
\left( \varkappa \right) _{\sigma }^{n}$ for all $n<\omega $ and all $\sigma
<\varkappa $. Any weakly compact cardinal $\varkappa $ is larger then the
first inaccessible one, moreover, $\varkappa $ is inaccessible itself and is
the $\varkappa $'s inaccessible cardinal.

Hence, $\varkappa \rightarrow \left( \varkappa \right) _{\gimel _{\alpha
}}^{n}$ for each $n<\omega $ and $\alpha <\iota $, where $\iota $ is the
first inaccessible cardinal.

\begin{theorem}
Let $\varkappa $ be a weakly compact cardinal; $X\in \mathcal{B}$ and $%
dim(X)=\varkappa $. Then $X$ contains an unconditional $\varkappa $-sequence.
\end{theorem}

\begin{proof}
Let $\{x_{\gamma }:\gamma <\varkappa \}$ be a sequence of linearly
independent elements of $X$ that spans $X$. Fix $n<\omega $. Let $%
\{a_{0},...,a_{n-1}\}$ are scalars, $b\in \mathbb{R}^{+}$ - a positive real
number. There exists a continuum number of choosing of $\phi
=\{n,a_{0},...,a_{n-1}\}$. The set $\Phi $ of such elements $\phi $ may
assumed to be well-ordered; let $\alpha _{\phi }$ be the ordinal number of $%
\phi $ in this ordering.

It will be said that a finite sequence $\gamma _{0}<\gamma _{1}<...<\gamma
_{n}<\varkappa $ belongs to a class $C_{\alpha }^{b}$ if
\begin{equation*}
\left\| \sum_{k=0}^{n}a_{k}x_{\gamma _{k}}\right\| =b.
\end{equation*}

For a fixed $\alpha $ a set of all classes $C_{\alpha }^{b}$ forms a
partition of $[\varkappa ]^{n}$ on continuum parts. Since $\varkappa $ is
weakly compact, there exists a homogeneous subset $K(\alpha )\subset
\varkappa $ for this partition; $card(K(\alpha ))=\varkappa $.

Next we proceed by the induction. Assume that sets $\{K(\beta ):\beta
<\alpha \}$ are already chosen. If $\alpha =\gamma +1$ for some $\gamma $
then $K(\alpha )$ is defined to be a homogeneous subset of $K(\gamma )$ with
respect to the partition $\{C_{\gamma +1}^{b}:b\in \mathbb{R}^{+}\}$. If $%
\alpha $ is a limit ordinal then $K(\alpha )=\cap \{K(\beta ):\beta <\alpha
\}$. Each $K(\alpha )$ is of cardinality $\varkappa $. Clearly, $\cap
\{K(\beta ):\beta <\varkappa \}$ is of cardinality $\varkappa $ too, is well
ordered and, hence, is a desired spreading $\varkappa $-sequence.
\end{proof}

The converse result is also true.

\begin{theorem}
For every cardinal $\varkappa $, which is not weakly compact, there exists a
Banach space $X$ of dimension $dim(X)=\varkappa $, which does not contain
any unconditional $\varkappa $-sequence.
\end{theorem}

\begin{proof}
Let $W$ be an abstract set of cardinality $\varkappa $. On $W$ may be
defined a linear order, say, $\ll $, such that $\left\langle W,\ll
\right\rangle $ contains neither strictly increasing nor strictly decreasing
$\varkappa $-sequences with respect to this order. Let $X_{0}$ be a vector
space of all formal finite linear combinations of elements of $W$ with real
scalars. Elements of $X_{0}$ are formal sums $x=\sum_{i\in I}a_{i}w_{i}$,
where $I$ is finite. Let $\{y_{n}:n<\omega \}$ be a spreading sequence of
elements of some Banach space $Y$. It will be assumed that $\{y_{n}:n<\omega
\}$ is not equivalent to any symmetric sequence.

For any element $x=\sum_{i\in I}a_{i}w_{i}$ let $\left| \left\| x\right\|
\right| =\left\| \sum_{i\in I}a_{i}z_{i}\right\| _{Y}$. Let $X$ be a
completition of $X_{0}$ with respect to this norm. Of course, $%
dim(X)=\varkappa $. Clearly, $\{w_{w}:w\in W\}$ is a spreading family,
indexed by itself. If $\{s_{\alpha }:\alpha <\varkappa \}$ is a spreading $%
\varkappa $-sequence of elements of $X$, then it must be equivalent to a
some block-sequence with respect to $\{w_{w}:w\in W\}$. However this is
impossible for our choosing of $W$.
\end{proof}

This result may be improved. However, a power of the improvement essentially
depends on a model of the set theory, which lies in a base of the whole
Banach space theory.

\begin{theorem}
For any cardinal $\tau \geq \omega $ there exists a Banach space $X$ of
dimension $dim(X)=2^{\tau }$ (=$\gimel _{1}\left( \tau \right) $) that does
not contain any subsymmetric $\tau ^{+}$-sequences.
\end{theorem}

\begin{proof}
A ''Banach'' part of the theorem is similar to the previous one. A
set-theoretical part is to find on a set $I$ of a cardinality $2^{\tau }$ a
such linear order, say $\ll $, that $\left\langle I,\ll \right\rangle $
contains neither strictly increasing nor strictly decreasing $\tau ^{+}$%
-sequences. If $\tau =\omega $, for $\ll $ may be chosen the natural order
on the set $\mathbb{R}$ of real numbers; if $\tau >\omega $, as $\ll $ may
be regarded a lexicographical ordering on a set $^{\tau }2$ (of all
functions from $\tau $ to a set $\{0,1\}$).
\end{proof}

\section{A classification of Banach spaces}

Preceding results shows that it is of sense to introduce a partition of the
class $\mathcal{B}$ of all Banach spaces by sets
\begin{equation*}
\mathcal{B}_{\varkappa }=\{Y\in \mathcal{B}:dim(Y)=\varkappa \},
\end{equation*}
which may be regarded as an analogue of horizontal strips of some
generalized coordinate system for Banach spaces. Below were presented
results, which demonstrate dependence of properties of Banach spaces on
their dimension.

As vertical strips of this generalized coordinate system may be regarded
classes of finite equivalence.

\begin{definition}
Let $X$, $Y$ be Banach spaces. $X$ is said to be \textit{finitely
representable} in $Y$, shortly: $X<_{f}Y$, if for every $\varepsilon >0$ and
for every finite dimensional subspace $A$ of $X$ there is a subspace $B$ of $%
Y$ and an isomorphism $u:A\rightarrow B$ such that $\left\| u\right\|
\left\| u^{-1}\right\| <1+\varepsilon $. $X$ and $Y$ are said to be\textit{\
finitely equivalent}, $X\sim _{f}Y$, if $X<_{f}Y$ and $Y<_{f}X$. So, every
Banach space $X$ generates a \textit{class of finite equivalence}
\begin{equation*}
X^{f}=\{Y\in \mathcal{B}:X\sim _{f}Y\}.
\end{equation*}
\end{definition}

It will be shown that for any infinite dimensional Banach space $X$ the
corresponding class $X^{f}$ is proper: it contains spaces of arbitrary large
dimension.

For a proof there will be used ultrapower of Banach spaces.

\begin{definition}
Let $I$ be a set; $Pow(I)$ be a set of all its subsets. An ultrafilter $D$
over $I$ is a subset of $Pow(I)$ with following properties:

\begin{itemize}
\item  $I\in D$;

\item  If $A\in D$ and $A\subset B\subset I$, then $B\in I$;

\item  If $A$, $B$ $\in D$ then $A\cap B\in D$;

\item  If $A\in D$, then $I\backslash A\notin D$.
\end{itemize}
\end{definition}

\begin{definition}
Let $I$ be a set; $D$ be an ultrafilter over $I$; $\{X_{i}:i\in I\}$ be a
family of Banach spaces. An \textit{ultraproduct } $(X_{i})_{D}$ is given by
a quotient space
\begin{equation*}
(X)_{D}=l_{\infty }\left( X_{i},I\right) /N\left( X_{i},D\right) ,
\end{equation*}
where $l_{\infty }\left( X_{i},I\right) $ is a Banach space of all families $%
\frak{x}=\{x_{i}\in X_{i}:i\in I\}$, for which
\begin{equation*}
\left\| \frak{x}\right\| =\sup \{\left\| x_{i}\right\| _{X_{i}}:i\in
I\}<\infty ;
\end{equation*}
$N\left( X_{i},D\right) $ is a subspace of $l_{\infty }\left( X_{i},I\right)
$, which consists of such $\frak{x}$'s$\ $that
\begin{equation*}
\lim_{D}\left\| x_{i}\right\| _{X_{i}}=0.
\end{equation*}
\end{definition}

If all $X_{i}$'s are all equal to a space $X\in \mathcal{B}$ then an
ultraproduct is said to be an \textit{ultrapower} and is denoted by $\left(
X\right) _{D}$.

An operator $d_{X}:X\rightarrow \left( X\right) _{D}$ that asserts to any $%
x\in X$ an element $\left( x\right) _{D}\in \left( X\right) _{D}$, which is
generated by a stationary family $\{x_{i}=x:i\in I\}$, is called the \textit{%
canonical embedding }of $X$ into its ultrapower $\left( X\right) _{D}$.

It is well-known that a Banach space $X$ is finitely representable in a
Banach space $Y$ if and only if there exists such ultrafilter $D$ (over $%
I=\cup D$) that $X$ is isometric to a subspace of the ultrapower $(Y)_{D}$
(cf. [9], an excellent exposition on ultrapowers in the Banach space theory.
Below will be presented some important results on ultrapowers that does not
contained in [9]).

\begin{definition}
Let $D$ be an ultrafilter, $\varkappa $ be an cardinal. An ultrafilter $D$
is said to: be

\begin{itemize}
\item  $\varkappa $-regular if there exists a subset $G\subset D$ of
cardinality $card(G)=\varkappa $ such that any $i\in I$ belongs only to a
finite number of sets $e\in G$;

\item  $\varkappa $-complete if an intersection of an nonempty subset $%
G\subset D$ of cardinality $<\varkappa $ belongs to $D$, i.e., if
\begin{equation*}
G\subset D\text{ \ and \ }card(G)<\varkappa \text{ \ implies that }\cap G\in
D;
\end{equation*}

\item  Principal (or non-free), if there exists $i\in I$ such that $%
D=\{e\subset I:i\in e\}$ (such ultrafilter is said to be generated by $i\in
I $);

\item  Free, if it is not a principal one;

\item  Countably incomplete if it is not $\omega _{1}$-complete.
\end{itemize}
\end{definition}

\begin{theorem}
Let $X$ be a Banach space, $dim(X)=\varkappa $; $D$ be an ultrafilter. The
canonical embedding $d_{X}:X\rightarrow \left( X\right) _{D}$ maps $X$
\textbf{on }$\left( X\right) _{D}$ if and only if $D$ is $\varkappa ^{+}$%
-complete.
\end{theorem}

\begin{proof}
Let $D$ be $\varkappa ^{+}$-complete. If $card(I)\leq \varkappa $ then $D$
is non-free and, hence, $d_{X}X=\left( X\right) _{D}$. Assume that $%
card(I)>\varkappa $ and that $\left( x_{i}\right) _{D}\in \left( X\right)
_{D}\backslash d_{X}X$. Of course, $\left( x_{i}\right) _{D}$ is generated
by a such family $\left( x_{i}\right) _{i\in I}$ that for any $\varepsilon
>0 $
\begin{equation*}
I_{0}(\varepsilon )=\{i\in I:\left\| x_{i}-x\right\| >\varepsilon \}\in D
\end{equation*}
for any $x\in X$. Consider a function $F:I_{0}\rightarrow X$, which is given
by $f(i)=x_{i}$. Since $dim(X)=\varkappa $, a partition
\begin{equation*}
I=\cup \{f^{-1}\left( x\right) :x\in I_{0}\left( \varepsilon \right) \}\cup
\{I\backslash I_{0}\left( \varepsilon \right) \}
\end{equation*}
participate $I$ on $\tau $ sets where $\tau <\varkappa ^{+}$. Since $D$ is $%
\varkappa ^{+}$-complete, one of sets of the partition belongs to $D$.
Because of $I_{0}\left( \varepsilon \right) \in D$, $I\backslash I_{0}\left(
\varepsilon \right) \notin D$. Hence there exists such $x\in X$ that $%
f^{-1}\left( x\right) \in D$. Clearly, this contradicts with the assumption $%
\left( x_{i}\right) _{D}\in \left( X\right) _{D}\backslash d_{X}X$.

Conversely, assume that $d_{X}X=\left( X\right) _{D}$. Since $%
dim(X)=\varkappa $ there exists $\varepsilon >0$ and a set $A\subset X$ such
that $card(A)=\varkappa $ and $\left\| a-b\right\| \geq \varepsilon $ for
all $a\neq b\in A$.

Let $I=\cup \{I_{a}:a\in B\}$ be a partition of $I$ on $\tau
=card(B)<\varkappa ^{+}$ parts. Let a function $F:I\rightarrow A$ is given
by:
\begin{equation*}
F(i)=a\text{ \ if and only if }i\in I_{a}\text{.}
\end{equation*}
Then $\left( F\left( i\right) \right) _{D}\in \left( X\right) _{D}=d_{X}X$,
and, hence, $\left( F\left( i\right) \right) _{D}=d_{X}\left( a\right) $ for
some $a\in A$. Because of $\left\| a-b\right\| \geq \varepsilon $ for all $%
a\neq b\in A$, $F^{-1}(a)\in D$. However, by our assumption, $%
F^{-1}(a)=I_{a} $. Thus, $I_{a}\in D$ and, since a partition of $I$ by $\tau
$ parts was arbitrary, $D$ is $\varkappa ^{+}$-complete.
\end{proof}

\begin{corollary}
Let $X$ be a Banach space, $D$ be an ultrafilter. then

\begin{enumerate}
\item  If $dim(X)<\omega $ then $\left( X\right) _{D}=X$;

\item  If $\omega \leq dim(X)<\upsilon $ where $\upsilon $ is the first
measurable cardinal, and $D$ is a free ultrafilter then $\left( X\right)
_{D}\backslash d_{X}X\neq \varnothing $;

\item  If $X$ is of infinite dimension and ultrafilter $D$ is countably
incomplete, then $\left( X\right) _{D}\backslash d_{X}X\neq \varnothing $
too.
\end{enumerate}
\end{corollary}

\begin{proof}
Since any ultrafilter is $\omega $-complete (by definitions), all results
follows from the preceding theorem.
\end{proof}

Nevertheless, by choosing an ultrafilter, a difference $\left( X\right)
_{D}\backslash d_{X}X$ may be maiden arbitrary large.

\begin{theorem}
Let $X\in \mathcal{B}_{\varkappa }$; $D$ be a $\tau $-regular ultrafilter
over a set $I$of cardinality $card(I)=\tau $. Then
\begin{equation*}
dim((X)_{D})=(dim(X))^{\tau }=\varkappa ^{\tau }.
\end{equation*}
\end{theorem}

\begin{proof}
Let $A^{\prime }\subset X$ be a set of cardinality $\varkappa $ which is
dense in $X$. Immediately,
\begin{equation*}
dim\left( \left( X\right) _{D}\right) \leq dim\left( l_{\infty }\left(
X,I\right) \right) \leq card\left( ^{I}A\right) =\left( card\left( A\right)
\right) ^{card\left( I\right) }=\varkappa ^{\tau }.
\end{equation*}

Assume now that $A\subset A^{\prime }$ is a set of the same cardinality $%
\varkappa $ such that for some $\varepsilon >0$ and any $a,b\in A$, $a\neq b$%
, $\left\| a-b\right\| \geq \varepsilon $. Let $\{a_{i}:i<\varkappa \}$ be a
numeration of elements of $A$. Any finite subset $s=\{x_{0},...,x_{n-1}\}%
\subset A$ spans a finite dimensional subspace $X_{s}$ of $X$. A set $B$ of
all subspaces of kind $X_{s}$ of $X$ is also of cardinality $\varkappa $.
Let $\{X_{\alpha }:\alpha <\varkappa \}$ be a numeration of $B$. Consider an
ultraproduct $\left( X_{\alpha }\right) _{D}$. Since $X_{\alpha
}\hookrightarrow X$ for all $\alpha <\varkappa $, it may be assumed that $%
\left( X_{\alpha }\right) _{D}\hookrightarrow \left( X\right) _{D}$ and,
hence, $dim\left( \left( X_{\alpha }\right) _{D}\right) \leq dim\left(
\left( X\right) _{D}\right) $. Thus, to prove the theorem it is enough to
show that $dim\left( \left( X_{\alpha }\right) _{D}\right) \geq \varkappa
^{\tau }$.

From $\tau $-regularity of $D$ follows that there exists a subset $G\subset
D $, $card(G)=\tau $, such that any $i\in I$ belongs only to a finite
numbers of sets $e\in D$. Define on $G$ a linear order, say, $\ll $. Let $%
g:G\rightarrow A$ be a function. Let a function $f_{g}:I\rightarrow B$ be
given by
\begin{equation*}
f_{g}\left( i\right) =span\{a_{g\left( e_{k}\left( i\right) \right)
}:k<n_{i}\}=B_{i},
\end{equation*}
where $n_{i}$ is a cardinality of a set of those sets $e\in G$, to which $i$
belongs; $\{e_{0}\left( i\right) \ll ...\ll e_{n_{i}}\left( i\right) \}$ is
a list of them. Let $h:G\rightarrow A$ be any other function: $h\neq g$.
Then there exists $e\in G$ such that $h\left( e\right) \neq g\left( e\right)
$. Corresponding functions $f_{g}$ and $f_{h}$ are differ also. Indeed, for
any $i\in e$ the set $e$ is contained in a finite sequence $e_{0}\left(
i\right) \ll ...\ll e_{n_{i}}\left( i\right) $ of all sets that contain $i$.
if its number in a such sequence if equal to $k$, then
\begin{equation*}
f_{g}\left( i\right) =span\{...,g\left( e_{k}\right) ,...\}\neq
span\{...,h\left( e_{k}\right) ,...\}=f_{h}\left( i\right) .
\end{equation*}

Note that $e\in D$ and that $f_{g}\left( i\right) \neq f_{h}\left( i\right) $
for all $i\in e$.

Hence $f_{g}\left( i\right) $ and $f_{h}\left( i\right) $ generate different
elements $\frak{x}\left( f\right) $ and $\frak{x}\left( h\right) $ of $%
\left( X_{\alpha }\right) _{D}$; moreover $\left\| \frak{x}\left( f\right) -%
\frak{x}\left( h\right) \right\| \geq \varepsilon $.

Since $card\left( ^{G}A\right) =\varkappa ^{\tau }$, and different elements $%
f,h\in ^{G}A$ generate different elements $\frak{x}\left( f\right) $ and $%
\frak{x}\left( h\right) $ of $\left( X_{\alpha }\right) _{D}$ with $\left\|
\frak{x}\left( f\right) -\frak{x}\left( h\right) \right\| \geq \varepsilon $%
, it is clear that $dim\left( \left( X_{\alpha }\right) _{D}\right) \geq
\varkappa ^{\tau }$.
\end{proof}

\begin{remark}
It may be proved that for any infinite-dimensional Banach space $X$ and any
countably incomplete ultrafilter $D$,
\begin{equation*}
dim((X)_{D})=(dim((X)_{D}))^{\omega }.
\end{equation*}
Hence, by using ultrapowers, cannot be obtained any space $(X)_{D}$, which
dimension $\varkappa $ is of countable confinality (i.e. such that $%
cf(\varkappa )=\omega $; e.g., $\aleph _{\omega }$; $\gimel _{\omega
_{\omega }}$ and so on).
\end{remark}

\begin{theorem}
For any infinite dimensional Banach space $X$ the corresponding class $X^{f}$
contains spaces of every infinite cardinality (i.e., for any cardinal $%
\varkappa \geq \omega $, $X^{f}\cap \mathcal{B}_{\varkappa }\neq \varnothing
$.
\end{theorem}

\begin{proof}
Let $dim(X)=\tau $. If $\varkappa =\omega $ then we choose a countable
sequence $\{A_{i}:i<\omega \}$ of finite dimensional subspaces of $X$ which
is dense in a set $H(X)$ of all different finite dimensional subspace of $X$
(isometric subspaces in $H(X)$are identified), equipped with a metric
topology, which is induced by the Banach-Mazur distance. Clearly, $%
X_{0}=span\{A_{i}:i<\omega \}\hookrightarrow X$ is separable and is finite
equivalent to $X.$

If $\varkappa \leq \tau $ then as a representative of $X^{f}\cap \mathcal{B}%
_{\varkappa }$ may be chosen any subspace of $X$ of dimension $\tau \ $that\
contains a subspace $X_{0}$.

If $\tau <\varkappa $ then, by the theorem 8 it may be found an ultrapower $%
(X)_{D}$ of dimension $\geq \varkappa $. Certainly, $(X)_{D}\sim _{f}X$.
Choose any subspace of $(X)_{D}$ of dimension $\varkappa $, which contains a
subspace $d_{X}X$.
\end{proof}

\section{Classes of $\protect\omega $-equivalence}

If one replace in the definition of finite representability finite
dimensional subspaces with infinite dimensional ones, it will be obtained a
new relation of $\omega $\textit{-representability}.

\begin{definition}
Let $X$, $Y$ be Banach spaces. $X$ is said to be $\omega $-representable in $%
Y$, shortly: $X<_{\omega }Y$, if for any separable subspace $%
E\hookrightarrow X$ there exists a separable $F\hookrightarrow Y$, which is
isomorphic to $E$. $X$ and $Y$ are $\omega $-equivalent, $X\sim _{\omega }Y$%
, if $X<_{\omega }Y$ and $Y<_{\omega }X$. Any $X\in \mathcal{B}$ generates a
class $X^{\omega }$ of all Banach spaces that are $\omega $-equivalent to $X$%
:
\begin{equation*}
X^{\omega }=\{Y\in \mathcal{B}:X\sim _{\omega }Y\}
\end{equation*}
\end{definition}

The theorem 9 is not remain valid if one replaces $X^{f}$ with $X^{\omega }$%
. Of course, for some Banach spaces $X$ a class $\ X^{\omega }$ also may
contain spaces of arbitrary large dimension (e.g. if $X=l_{p}$ or $X=L_{p}$
for some $p\in \lbrack 1,\infty ]$).

It is easy to show that any Banach space $X$, which contains a subspace with
a symmetric basis, generates a class $X^{\omega }$, which contains spaces of
arbitrary infinite dimension.

However, a method of ultrapowers cannot be used in a general case to
construct other then $X$ spaces from a class $X^{\omega }$. Indeed, it is
known (cf. [11]) that for any infinite dimensional Banach space $X$ its
ultrapower $\left( X\right) _{D}$ by a countably incomplete ultrafilter $D$
is universal for all separable Banach spaces $Y$ which are finitely
representable in $X$. Hence, the following proposition is true.

\begin{proposition}
For any infinite dimensional Banach space and each pair of countably
incomplete ultrafilters $D$ and $E$ spaces $\left( X\right) _{D}$ and $%
\left( X\right) _{E}$ are $\omega $-equivalent.
\end{proposition}

Another example gives a hereditarily indecomposable reflexive separable
space $X_{HI}$, that was constructed in [2]. Obviously, a class $\left(
X_{HI}\right) ^{\omega }$ cannot contain any non-separable space (since any
non-separable reflexive Banach space is decomposable).

\begin{remark}
Let $X$ be a Banach space of dimension $\varkappa \geq \upsilon $, where $%
\upsilon $ is the first measurable cardinal, and $D$ be an $\omega $%
-complete free ultrafilter It may be proved that every subspace of $\left(
X\right) _{D}$ of dimension, strictly less then $\upsilon $, is isometric to
a subspace of $X$. This means that in a such case $X\sim _{\omega }\left(
X\right) _{D}$.

Since a dimension of $\left( X\right) _{D}$ (by choosing of an $\omega $%
-complete free ultrafilter) may be maiden arbitrary large, it may be
conclude that if a dimension $dim\left( X\right) \geq \upsilon $, where $%
\upsilon $ is the first measurable cardinal, then a class $X^{\omega }$
contains spaces of any dimension $\geq \gimel _{1}$.\
\end{remark}

The (too large) dimension $\upsilon $ may be essentially lowered.

\begin{theorem}
For any Banach space $X$ of dimension $dim(X)\geq \gimel _{\omega _{1}}$ and
any cardinal $\varkappa \geq \gimel _{\omega _{1}}$, the corresponding class
$X^{\omega }$ contains a space $Y_{\varkappa }$ of dimension $\varkappa $.
\end{theorem}

The proof of this theorem, which involves some model theory (and a series of
additional definitions) is contained in an appendix to this article..

\begin{remark}
From this theorem follows immediately that if a cardinal $\varkappa \left(
\omega \right) $ exists then any Banach space $X$ of dimension $%
dim(X)=\gimel _{\omega _{1}}$ contains a subspace with an unconditional
basis. Indeed, since $\gimel _{\omega _{1}}<\varkappa \left( \omega \right) $%
, from the equality $dim(X)=\gimel _{\omega _{1}}$ \ follows that the class $%
X^{\omega }$ contains a space $Y$ of dimension $\varkappa \left( \omega
\right) $. This space contains a subspace $W$ with an unconditional basis
(see [4]). Since $Y<_{\omega }X$, $X$ also has a subspace isomorphic to $W$.
\end{remark}

\begin{remark}
Note that Theorem 5 makes more clear the role of the cardinal $\gimel
_{\omega _{1}}$.
\end{remark}

\section{Appendix}

\subsection{Some model theory.}

To express the notion of a Banach space in the logic consider a \textit{%
language }$\frak{L}$ - a set of symbols, which includes besides of logical
symbols (prepositions $\&,\vee ,\urcorner $, quantifiers $\forall ,\exists $
and variables $u,v,x_{1},x_{2},\ $etc) also non-logical primary symbols:

\begin{itemize}
\item  A binary functional symbol $+$;

\item  A countable number of unary functional symbols $\{f_{q}:q\in \mathbb{Q%
}\}$;

\item  An unary predicate symbol $B$.
\end{itemize}

Any Banach space $\frak{X}=\left\langle X,\left\| \cdot \right\|
\right\rangle $ may be regarded as a model for all logical propositions (or
formulae) that are satisfied in $X$. A set $\left| \frak{X}\right| =X$ is
called a support (or an absolute) of a model $\frak{X}$, in which all non
logical symbols of the language $\frak{L}$ are interpreted as follows:

\begin{itemize}
\item  $+^{\frak{X}}$ is interpreted as the addition of vectors from $X$;

\item  ($f_{q})^{\frak{X}}$ for a given $q\in \mathbb{Q}$ is interpreted as
a multiplication of vectors from $X$ by a rational scalar $q$

\item  $B^{\frak{X}}$ is interpreted as the unit ball $B(X)=\{x\in X:\left\|
x\right\| \leq 1\}$.
\end{itemize}

The language $\frak{L}$ was introduced by J. Stern [11], who used for
examine the Banach space\ theory the \textit{first order logic}. In this
logic any quantifier acts only on elements of a supports of corresponding
structures (e.g., a quantification over either all countable subsets of $X$
or over all natural $n$'s is forbidden); only formulae that contains only
finite strings of symbols are admissible. We shall used a less restrictive%
\textit{\ infinitary logic }$\frak{L}_{\omega _{1},\omega }$, which may be
obtained from the first order logic by the addition of two new rules.

\begin{enumerate}
\item  \textit{If} $\Phi $ \textit{is a set of formulae of }$\frak{L}%
_{\omega _{1},\omega }$, \textit{such that} $card\left( \Phi \right) \leq
\omega $,\textit{\ then }$\&\Phi $ \textit{also is a formula of }$\frak{L}%
_{\omega _{1},\omega }$;

\item  \textit{If }$\varphi $ \textit{is a formula of} $\frak{L}_{\omega
_{1},\omega }$ \textit{and }$V$ \textit{is a finite set of variables then }$%
\left( \forall V\right) \varphi $ \textit{is a formula of} $\frak{L}_{\omega
_{1},\omega }$.
\end{enumerate}

In the logic $\frak{L}_{\omega _{1},\omega }$ all axioms of Banach spaces
can be expressed instead of the axiom of completeness (its formulation
requires a quantification over a countable set of variables). For this
reason (and for advantages of a countable language) any Banach space $X$
will be regarded as a $\mathbb{Q}$-vector normed space. Obviously, to any
pair of $\mathbb{Q}$-vector normed spaces $\frak{X}$, $\frak{Y}$, which are
isomorphic in a model theoretic sense (i.e. that are such that there exists
an one-to-one map $u:\left| \frak{X}\right| \rightarrow \left| \frak{Y}%
\right| $, which holds all functions and predicates of $\frak{L}$)
corresponds (under their completition, if it needed) a pair of isometric
Banach spaces. Thus, the infinitary logic $\frak{L}_{\omega _{1},\omega }$
will be sufficient for our needs.

If $\varphi $ is a formula of $\frak{L}_{\omega _{1},\omega }$, which is
satisfied in a Banach space $X$, we shall write $X\vDash \varphi $. Let
\begin{equation*}
Th(X)=\{\varphi \in \frak{L}_{\omega _{1},\omega }:X\vDash \varphi \}.
\end{equation*}

Banach spaces $X$ and $Y$ are said to be \textit{elementary equivalent in
the logic} $\frak{L}_{\omega _{1},\omega }$, shortly, $X\equiv _{\omega
_{1},\omega }Y$, if $Th(X)=Th(Y)$.

The classical Scott's theorem asserts that for any countable system$\frak{\ M%
}$ \ in a countable language for the logic $\frak{L}_{\omega _{1},\omega }$
there exists such $\frak{L}_{\omega _{1},\omega }$-formula $\psi \left(
\frak{M}\right) $ that the condition $\frak{N}\vDash \psi \left( \frak{M}%
\right) $ for any other system $\frak{N}$ \ implies that systems $\frak{M}$
and$\ \frak{N}$ are isomorphic. In our case, i.e., in the Banach space
setting, this theorem means that for any separable Banach space $X$ the
condition $\ X\equiv _{\omega _{1},\omega }Y$ implies that $X$ and $Y$ are
isometric.

As was noted by O.V. Belegradek (a personal communication), only small
changes in the proof of the Scott's theorem are needed to do, to obtain the
following its improvement:

\textit{For any system}$\frak{\ M}$ \textit{in a countable language for the
logic} $\frak{L}_{\omega _{1},\omega }$ \textit{there exists} \textit{such }$%
\frak{L}_{\omega _{1},\omega }$\textit{-formula }$\xi \left( \frak{M}\right)
$\textit{\ that a condition }$\frak{N}\vDash \psi \left( \frak{M}\right) $%
\textit{\ for any other system }$\frak{N}$ \textit{implies that systems }$%
\frak{M}$ and$\ \frak{N}$\textit{\ have isomorphic countable subsystems.}

From this immediately follows the following result.

\begin{theorem}
Let $X$, $Y$ are Banach spaces. If $X\equiv _{\omega _{1},\omega }Y$\ then $%
X\sim _{\omega }Y$.
\end{theorem}

\begin{proof}
From the preceding improvement of the Scott's theorem follows that $X\sim
_{\omega }Y$ if and only if $X\vDash Th(Y)_{\forall }$ and $Y\vDash
Th(X)_{\forall }$, where for any theory $T$ (i.e. for any set of $\frak{L}%
_{\omega _{1},\omega }$-formulae) $T_{\forall }$ denotes a set of all its
universal consequences. Of course, $T_{\forall }\subseteq T$ for any theory $%
T$. Therefore $X\equiv _{\omega _{1},\omega }Y$ implies that $X\sim _{\omega
}Y$.
\end{proof}

The next notion that will be needed is a notion of \textit{a Hanf number }of
the logic $\frak{L}_{\omega _{1},\omega }$.

Let $X$ be a $\mathbb{Q}$-vector normed space. Let $ModTh(X)=\{Y:Y\equiv
_{\omega _{1},\omega }X\}$, where $Y$ is a $\mathbb{Q}$-vector normed space
too. The Hanf number of the logic $\frak{L}_{\omega _{1},\omega }$ is the
least cardinal $h\left( \frak{L}_{\omega _{1},\omega }\right) $ with the
property:

\textit{If a model} $\frak{M}$ \ \textit{for this logic is of cardinality} $%
card\left| \frak{M}\right| \geq h\left( \frak{L}_{\omega _{1},\omega
}\right) $, \textit{then for any} \textit{cardinal} $\varkappa \geq h\left(
\frak{L}_{\omega _{1},\omega }\right) $ t\textit{here exists a model} $\frak{%
N}\equiv _{\omega _{1},\omega }\frak{M}$ \textit{of cardinality} $card\left|
\frak{N}\right| =\varkappa $.

In our case this means that for each Banach space $X$ of dimension $%
dim(X)\geq h\left( \frak{L}_{\omega _{1},\omega }\right) $ and for any $%
\varkappa \geq h\left( \frak{L}_{\omega _{1},\omega }\right) $ there exists
a Banach space $Y\equiv _{\omega _{1},\omega }X$ of dimension $%
dim(Y)=\varkappa $.

The exact meaning of the Hanf number for the logic $\frak{L}_{\omega
_{1},\omega }$ is known (see [13]):
\begin{equation*}
h\left( \frak{L}_{\omega _{1},\omega }\right) =\gimel _{\omega _{1}}
\end{equation*}

Now we are ready to present a proof of the Theorem 7.

\subsection{Proof of the theorem 10.}

\begin{proof}
Let $X$ be a Banach space of dimension $dim(X)=\gimel _{\omega _{1}}$.
Consider $X$ to be a model for the infinitary logic $\frak{L}_{\omega
_{1},\omega }$. Since its dimension is equal to the Hanf number $h\left(
\frak{L}_{\omega _{1},\omega }\right) $, for any $\varkappa \geq \gimel
_{\omega _{1}}$, then there exists $Y\equiv _{\omega _{1},\omega }X$ of
dimension $dim(Y)=\varkappa $. Clearly, the condition $Y\equiv _{\omega
_{1},\omega }X$ implies that $Y\sim _{\omega }X$. The last equivalence
proves the theorem.
\end{proof}

\begin{remark}
Let $X$ be a Banach space of dimension $\varkappa >\omega _{1}$. Then there
exists a Banach space $Y\equiv _{\omega _{1},\omega }X$ of dimension $%
dim(Y)=\omega _{1}$. This follows from the general result of R.L.~Vaught
[13]: the L\"{o}venheim number $l\left( \frak{L}_{\omega _{1},\omega
}\right) $ for a countable language $\frak{L}$ ( i.e. the least cardinal $%
\tau $ such that every $\frak{L}$-model $\frak{A}$ of cardinality $\varsigma
>\tau $ has a subsystem $\frak{B}\vDash Th\left( \frak{A}\right) $ of
cardinality $\tau $) is equal to $\omega _{1}$.
\end{remark}

\begin{acknowledgement}
Author wish to express his thanks to O.V. Belegradek for helpful advises and
recommendations
\end{acknowledgement}

\section{References}

\begin{enumerate}
\item  Lindenstrauss J. \textit{The geometric theory of the classical Banach
spaces,} Actes du Congres Intern. Math. Nice 1970, v. \textbf{2 }(1971)
365-372

\item  Gowers W.T., Maurey B.\textit{\ The unconditional basic sequence
problem}, Journal AMS \textbf{6} (1993) 851-874

\item  Ketonen J. \textit{Banach spaces and large cardinals,} Fund Math.
\textbf{81 }(1974) 291-303

\item  Tokarev E.V. \textit{Notion of indiscernibless and on the p-Mazur
property in Banach spaces} (transl. from Russian) Ukrainian Mathematical
Journal \textbf{37} (1985) p.180-184

\item  Tokarev E.V. \textit{On one P. Terenzi's question}, Deposed in VINITI
no.\textbf{325-B88} (1988) 1-7 (in Russian)

\item  Terenzi P. \textit{On the Banach spaces of large density character},
Riv. Mat. Univ. Parma \textbf{12:4 }(1986) 111-120

\item  Erd\"{o}s P., Rado R. \textit{A partition calculus in set theory},
Bull. AMS \textbf{62} (1956) 427-489

\item  Brunel A., Sucheston L. \textit{On J - convexity and some ergodic
super properties of Banach spaces}, Trans. AMS \textbf{204} (1975) 79-90

\item  Heinrich S. \textit{Ultraproducts in Banach space theory}, J. Reine
Angew. Math. \textbf{313 }(1980) 72-104

\item  Chang C. C., Keisler H. J. \textit{Model Theory}, Amsterdam:
North-Holland 1973

\item  Stern J.\textit{\ The problem of envelopes for Banach spaces}, Israel
J. Math. \textbf{24:1} (1976) 1-15

\item  Hanf W. \textit{Incompactness in languages with infinitely long
expressions}, Fund. Math. \textbf{13} (1964) 309-324

\item  Vaught R.L. \textit{Denumerable models of complete theories}, in:
Infinistic Methods; Pergamon Press (1961) 303-321
\end{enumerate}

\end{document}